\documentclass{PMmodified}
\usepackage{amssymb, amsmath, amscd, latexsym, mathrsfs}
\input xy
\xyoption{all}
\address[r.reis@abdn.ac.uk]{Dept. of Mathematics, University of Aberdeen, AB24 3UE, United Kingdom}%
\address[m.weiss@abdn.ac.uk\newline\newline\newline\date{\today}]{Dept. of Mathematics, University of Aberdeen, AB24 3UE, United Kingdom }%

\newtheorem{thm}{Theorem}[section]
\newtheorem{cor}[thm]{Corollary}
\newtheorem{lem}[thm]{Lemma}

\newtheorem{prop}[thm]{Proposition}
\newtheorem{defn}[thm]{Definition}

\newtheorem{rem}[thm]{Remark}

\newtheorem{expl}[thm]{Example}

\numberwithin{equation}{section}

\newcommand{\lra}{\longrightarrow}
\newcommand{\co}{\colon\!}

\newcommand{\smin}{\smallsetminus}

\newcommand{\sha}{^{\sharp}}

\newcommand{\op}{^{\textup{op}}}
\newcommand{\nin}{\noindent}

\newcommand{\id}{\textup{id}}
\newcommand{\im}{\textup{im}}

\newcommand{\ad}{\textup{ad}}

\newcommand{\sM}{\mathcal M} 

\newcommand{\sE}{\mathcal E} 
\newcommand{\sF}{F} 
\newcommand{\spaJ}{P_*} 

\newcommand{\RR}{\mathbb R}

\received[]{}

\title[Maps to the plane I]{Smooth maps to the plane and Pontryagin classes \\ Part I: Local aspects}

\author[Rui Reis and Michael Weiss]{Rui Reis and Michael Weiss
\thanks{This project was generously supported by the Engineering and Physical Sciences Research Council (UK), Grant EP/E057128/1.}}

\begin{document}

\maketitle

\begin{abstract} We classify the most common local forms of smooth maps from a smooth manifold $L$ to the
plane. The word \emph{local} can refer to locations in the source $L$, but also to locations
in the target. The first point of view leads us to a classification of certain germs of maps, which we review here
although it is very well known. The second point of view leads us to a classification of certain \emph{multigerms} of maps.
\end{abstract}
\maketitle


\section{Introduction}
Our goal is to investigate locally uncomplicated smooth maps from a smooth
manifold $L$ of dimension $n+2$ to the plane $\RR^2$. Where we use the word \emph{local}, as in
\emph{locally uncomplicated}, we sometimes refer to locations in the source $L$,
sometimes to locations in the target $\RR^2$. The emphasis is on families of smooth maps; this is in contrast to
Morse theory, where the study of individual (locally uncomplicated) smooth maps from a manifold to
$\RR$ is a central topic. We are guided by two observations.

\medskip\nin
(i) Let $X$ be an open subspace of the space of all smooth maps $L\to \RR^2$ defined by prohibiting
certain singularities. It is a special case of a theorem due to Vassiliev \cite{Vassiliev1},\cite{Vassiliev2}
that $X$ has an accessible homotopy type or homology type if, loosely speaking, every smooth map $L\to \RR^2$ can be
approximated by a map which belongs to $X$, and moreover every smooth one-parameter family of smooth maps
$L\to \RR^2$ can be approximated by a path in $X$. Therefore we are inclined to define notions of locally uncomplicated
map $L\to \RR^2$ by prohibiting certain singularities or singularity types corresponding to a subset of an
appropriate jet space whose codimension in the jet space is at least $n+4$.

\medskip\nin
(ii) More restrictive notions of locally uncomplicated
map $L\to \RR^2$ can be obtained by prohibiting, for every $r\ge 1$,
certain configurations of $r$ singularities (\emph{multigerms}) in the source $L$, with the same image point in $\RR^2$.
The Vassiliev theorem mentioned above can be adapted to this setup \cite{RWIII}, although it is considerably harder to say which multigerms can be prohibited without making the resulting space of locally uncomplicated smooth maps
$L\to \RR^2$ homologically or homotopically inaccessible.

\medskip
These two observations raise two elementary classification problems, one for uncomplicated germs and one for
uncomplicated multigerms. The solution of the first problem is well known, but we review it. In the second problem, it is not straightforward to come up with a manageable interpretation of \emph{classification}. We propose one and
describe our solution.

\section{Germs of maps from the plane to the plane}\label{sec-lowleftright}

The classification of the most common map germs from plane to plane
up to left-right equivalence is well known. See
for example \cite{AicardiOhmoto}. (We are talking about smooth
map germs $f$ from $(\RR^2,0)$ to $(\RR^2,0)$. Two such germs $f_0$, $f_1$ are \emph{left-right
equivalent} if there exist diffeomorphism germs $\psi\co(\RR^2,0)\to (\RR^2,0)$ and
$\sigma\co(\RR^2,0)\to(\RR^2,0)$ such that $f_1=\sigma f_0\psi^{-1}$.)  We will repeat it here nevertheless and see some
normal forms and tell the story of each singularity type.

\subsection{Classification}
There are six types that we consider worthy of attention:
\emph{regu\-lar}, \emph{fold}, \emph{cusp}, \emph{swallowtail}, \emph{lips}
and \emph{beak-to-beak}. The regular (alias nonsingular) type is well understood.
The remaining five types are of rank 1, that is, the derivative at the origin has rank 1.
(The cases where the derivative has rank 0 are uninteresting to us because their codimension is
at least 4.) Among these, it is natural to distinguish between those for which the 1-jet prolongation is
transverse to the rank 1 stratum (fold, cusp and swallowtail) and those for which it is not
(lips and beak-to-beak). In the transverse case, the singularity set in the source is a
smooth curve in the plane, passing through the origin; in the non-transverse case, it is in some
way or other a singular curve, as we will see.

\smallskip\nin
{\it Fold}: The normal form is $f(x,y) = (x,y^2)$.
The singularity set in the source is a line (in the normal form, the $x$-axis) and the
singularity set in the target is also a line (in the normal form, again the $x$-axis).
The intrinsic second derivative \cite{GolubitskyGuillemin73}
at the origin is a nondegenerate quadratic form (defined
on the kernel of the first derivative, and with values in the cokernel of the first
derivative).

\smallskip\nin
{\it Cusp}: Normal form $f(x,y)=(x,y^3+xy)$.
The derivative matrix for the normal form is
\begin{align*}   df(x,y) & = \begin{bmatrix} 1 & 0 \\ y & 3y^2+x \end{bmatrix} \end{align*}
with determinant $(x,y)\mapsto 3y^2+x$. Hence the singularity set $\Sigma$ in the source is
the trajectory of $t\mapsto(-3t^2,t)$, a parabola. The singularity set in the target
is the trajectory of $t\mapsto (-3t^2,-2t^3)$.

\smallskip\nin
{\it Swallowtail}: Normal form $f(x,y)=(x,y^4+xy)$.
The singularity set $\Sigma$ in the source is the trajectory of $t\mapsto(-4t^3,t)$.
The singularity set in the target is the trajectory of $t\mapsto (-4t^3,-3t^4)$.\smallskip

\smallskip\nin
{\it Lips}: Normal form $f(x,y)=(x,y^3+x^2y)$. The singularity
set in the source is the set of zeros of the quadratic form $(x,y)\mapsto x^2+3y^2$, that is,
a single point. It is a manifold but it does not have dimension 1.\smallskip

\smallskip\nin
{\it Beak-to-beak}: Normal form $f(x,y)=(x,y^3-x^2y)$.
The singularity set in the source is the set of zeros of the quadratic form $(x,y)\mapsto x^2-3y^2$,
that is, the union of the lines described by $x=cy$ and $x=-cy$, where $c=3^{1/2}$. It has dimension 1
but it is not a manifold. The singularity set in the target is the union of the trajectories of
\[  t\mapsto (ct,-2t^3),\qquad t\mapsto (-ct,-2t^3). \]

\begin{rem} {\rm In all these formulae, the first coordinate $f_1$ of $f$ is $(x,y)\mapsto x$.
The best way to understand the classification is to regard the second coordinate $f_2$ of $f$
as an \emph{unfolding} of a germ $g:(\RR,0)\to (\RR,0)$, with unfolding parameter $x$.
The formula for $g$ can be seen by setting $x=0$. This gives $g(y)=y$ for the regular case,
$g(y)=y^2$ for the fold, $g(y)=y^3$ for cusp, lips and beak-to-beak, and $g(y)=y^4$ for the
swallowtail.\footnote{Catastrophe theory has names for these germs $g$ which sometimes clash with our names
for the corresponding maps $f$. The catastrophe theory names tend to describe the projection from the fiberwise
singularity set of the miniversal unfolding of $g$ to the parameter space of the unfolding.}
Each of the unfoldings can be pulled back from a miniversal
unfolding with parameter space $V$. The miniversal unfoldings are
as follows:
\begin{eqnarray}
g(y)=y^2: &  y^2 \\
g(y)=y^3: &  y^3+uy \\
\label{2.7.09:2} g(y)=y^4: & y^4-uy^2+vy~.
\end{eqnarray}
This is essentially in the notation of \cite[ch.15]{Broecker}, although we use $y$ where \cite{Broecker} has $x$.
The decisive features of the germs $f$ are therefore as follows:
\begin{itemize}
\item[(i)] the corresponding germ $g:(\RR,0)\to (\RR,0)$ obtained by setting $x=0$ in the formula for $f_2$~;
\item[(ii)] the smooth map $e:(\RR,0)\to (V,0)$ (where $V$ parametrizes the miniversal unfolding of the appropriate $g$)
such that $f_2$ as an unfolding is isomorphic to $e^*$ of the miniversal unfolding. This $e$ is in most
cases far from unique.
\end{itemize}
For us, $V=\RR$ or $V=\RR^2$. In the notation of \cite[ch.15]{Broecker}, the maps $e$ are as follows:
$e(x)=x\in\RR$ for the cusp, $e(x)=(0,x)\in\RR^2$ for the swallowtail, $e(x)=x^2\in\RR$ for the lips and
$e(x)\mapsto -x^2\in \RR$ for beak-to-beak.
}
\end{rem}

\medskip
It is not completely trivial to justify this classification. What the above arguments prove beyond doubt
is that we have a \emph{surjective} map from isomorphism classes of 1-parameter unfoldings
of germs $g$ (such as $g(y)=y^n$, with $n=1,2,3,4$) to the set of left-right equivalence classes
of germs $f:(\RR^2,0)\to (\RR^2,0)$ whose derivative at $0$ has rank 1. What remains to be done
is roughly the following:
\begin{itemize}
\item[(i)] to produce a ``sufficiently big'' list of
some of the $1$-parameter unfoldings of the germs $g$,
and to determine the corresponding
left-right equivalence classes of map germs $f\co(\RR^2,0)\to (\RR^2,0)$~;
\item[(ii)] to show that each of these left-right equivalence classes has codimension $\le 3$ and that
all remaining germs $(\RR^2,0)\to (\RR^2,0)$ taken together make up a subset of codimension $\ge4$.
\end{itemize}

\bigskip\nin
\subsection{Unfoldings}
We start with the list of unfoldings.
Every 1-parameter unfolding of a smooth function germ
$g:(\RR,0)\to(\RR,0)$ with nonzero Taylor series is isomorphic
to $e^*$ of the miniversal unfolding, where
\[ e:(\RR,0)\to (V,0) \]
is smooth and $V=V_g$ is the parameter space
for the miniversal unfolding of $g$. The fact that $e$ is usually not unique makes the classification difficult.
However, some special cases are easy.

\smallskip\nin
If $g(y)=y^2$, then $V_g$ is zero-dimensional.

\smallskip\nin
If $g(y)=y^3$, then $V_g$ is 1-dimensional. The proposed normal forms for $e$ are
$e(x)=x$, $e(x)=x^2$ and $e(x)=-x^2$. If $q:(\RR,0)\to \RR=V_g$ has
nonzero first derivative, then we can find an invertible $h:(\RR,0)\to(\RR,0)$ such that
$q=eh$ where $e(x)=x$, and that can be used to produce the required isomorphism. Similarly, if $q$
has zero first derivative but strictly positive second derivative, then we can find an
invertible $h:(\RR,0)\to(\RR,0)$ such that $q=eh$ where $e(x)=x^2$.
Similarly, if $q$ has zero first derivative but strictly negative second derivative, then
we can find an invertible germ $h:(\RR,0)\to(\RR,0)$ such that $q=eh$ where $e(x)=-x^2$.

\smallskip\nin
So in fact the only difficult case is the case where $g(y)=y^4$. We use the miniversal
unfolding given by~(\ref{2.7.09:2}). Hence $V=V_g$ is 2-dimensional.
We want to focus on map germs $e:(\RR,0)\to (V,0)$ with nonzero first derivative,
not parallel to the $u$-axis. (The $u$-axis is a distinguished direction in $V$ because
it is parallel to the cusp in $V$ obtained by projecting the fiberwise singularity set
of the unfolding to $V$.) The corresponding 1-parameter
unfolding of $g(y)=y^4$ then has the form $y^4+e_1(x)y^2+e_2(x)y$
with $e_2'(0)\ne 0$. Using $e_2$ to transform the source of $e$,
we can reduce to a situation where $e_2(x)=x$. So we have
\[ (x,y) \mapsto y^4+p_xy^2+xy \]
where $p_x=e_1(x)$. From example~\ref{expl-unfoldswallowtail} we know that this
is left-right equivalent to $(x,y)\mapsto(y^4+xy)$, which is the swallowtail normal form.

\medskip
The rest of our classification task is easier. The five singularity types, represented by the five
normal forms above,
are easy to distinguish by geometric properties which are invariant under left-right equivalence.

\smallskip\nin
For the fold type, the singularity set $\Sigma$ in the source is a smooth submanifold
of dimension 1, and $f|\Sigma$ is an immersion (near 0).

\smallskip\nin
For the cusp and swallowtail, the singularity set $\Sigma$
in the source is still a smooth submanifold of dimension 1, but $f|\Sigma$ is not an immersion near $0$.
To distinguish cusp and swallowtail, it is enough to show that the curves
\[ t\mapsto (-3t^2,-2t^3),~~~~~~ t\mapsto (-4t^3,-3t^4) \]
are not left-right equivalent. This is obvious by looking at the second (intrinsic)
derivative~\cite{Broecker,GolubitskyGuillemin73} at the
origin, which is nonzero in the cusp case, zero in the swallowtail case.

\smallskip\nin
For the lips and beak-to-beak, the singularity set in the source is not a smooth submanifold of
dimension 1; it is a point in the lips case and a ``node'' (two crossing lines) in the beaks-to-beaks
case.

\subsection{Codimension and stratification}
We turn to the codimension and stratification analysis.
Among other things we want to determine the codimension of each of the six types described above,
and we want to show that all remaining singularity types taken together constitute a set of codimension $>3$.
We start by summarizing the analytic characterizations of the six types. We can always assume
that $f:(\RR^2,0)\to (\RR^2,0)$ has the form
\[ (x,y)\mapsto (x,f_2(x,y)) \]
and $\partial f_2/\partial x$
vanishes at $0$. In the singular case, we also assume that $\partial f_2/\partial y$ vanishes at $0$.
The following table describes the six types by means of conditions on the 4th Taylor polynomial of $f_2$.
The conditions typically state that some term in the Taylor polynomial has to be zero (z)
or nonzero (n). For example, the table states that in the case of
a cusp, the coefficients of $y$ and $y^2$ must be zero while the coefficients of $xy$ and $y^3$
must be nonzero (and there are no further conditions). \newline
In the ``other conditions'' column of the table,
$b_3$, $d_1$ and $d_2$ are the coefficients of $y^3$, $xy^2$ and $x^2y$ respectively.
The expression $3b_3d_2-d_1^2$ arises when we trade $xy^2$ terms for $x^2y$ terms, composing
with a diffeomorphism germ (in the source) of the form $(x,y)\mapsto(x,y-kx)$ for some constant $k$.

\[
\begin{array}{|c|c|c|c|c|l|l|}
\hline
\rule{0mm}{4mm} y & y^2 & y^3 & y^4 & xy & \textup{other conditions} & \textup{Name} \\ \hline
\rule{0mm}{4mm} n &     &     &     &    &         & \textup{regular} \\ \hline
\rule{0mm}{4mm} z & n   &     &     &    &         & \textup{fold}   \\ \hline
\rule{0mm}{4mm} z & z   & n   &     & n  &         & \textup{cusp}   \\ \hline
\rule{0mm}{4mm} z & z   & n   &     & z  & 3b_3d_2-d_1^2>0     & \textup{lips}  \\  \hline
\rule{0mm}{4mm} z & z   & n   &     & z  & 3b_3d_2-d_1^2<0     & \textup{beak-to-beak}  \\  \hline
\rule{0mm}{4mm} z & z   & z   & n   & n  &         & \textup{swallowtail} \\ \hline
\end{array}
\]

\medskip\noindent

\medskip
\begin{defn} {\rm Let $\spaJ$ be the  real vector space of polynomial
maps $\RR^2\to \RR^2$ (viewed as jets), of degree $\le 4$, with vanishing constant term. We write
\[  \spaJ=\spaJ^2\cup \spaJ^1\cup \spaJ^0 \]  
where $\spaJ^{i}$ consists of all those elements of $\spaJ$ whose
linear term has rank $i$. Let $W^{\spaJ}\subset \spaJ$
consist of the polynomials whose germ at the origin belongs to one of the types regular,
fold, cusp, swallowtail, lips or beak-to-beak. Thus
\[ \spaJ^{2} \subset W^{\spaJ} \subset \spaJ^{1}\cup \spaJ^{2}. \]
Let's also introduce $N\subset \spaJ^{1}$, the submanifold
of those $f$ which have the form $f(x,y)=(x,f_2(x,y))$ where $f_2$
has vanishing first derivative. \newline
For $\spaJ^{2}$ we also write $G$, because it is a Lie group.
The group $G$ acts on the left and right of $W^{\spaJ}$ by composition
of polynomial mappings (followed by truncation to degree $\le 4$).
In other words, $G\times G\op$ acts on $W^{\spaJ}$ by
$(\varphi,\psi)\cdot f=\varphi f\psi$, for $\varphi,\psi\in G$ and $f\in W^{\spaJ}$.
}
\end{defn}

Our classification attempts so far describe some orbits of this action
of $G\times G\op$ on $W^{\spaJ}$.
(In particular our classification of some germs $f:(\RR^2,0)\to
(\RR^2,0)$ up to left-right equivalence can be formulated in terms
of Taylor expansions at the origin, up to degree 4 at most.)
We now wish to show that $W^{\spaJ}$ is open, to determine the codimensions in $W^{\spaJ}$ of the
six orbits, and show that the complement of $W^{\spaJ}$ has codimension $\ge 4$ in $\spaJ$.
We have already convinced ourselves
that every $g\in \spaJ^{1}$ is left-right equivalent to some $f\in
N$. In other words, the restricted action map $G\times N\times
G\to \spaJ^{1}$ is onto. The following lemma makes this more
precise:

\begin{lem}
\label{lem-fiberbundle}
The restricted action map $G\times N\times G\to  \spaJ^{1}$ is a fiber bundle.
\end{lem}

\proof Let $E\subset G\times \spaJ^{1}$ be the smooth submanifold consisting of all pairs
$(\varphi,g)$, with $\varphi\in G$ and $g\in \spaJ^{1}$, such that the first derivative
of $\varphi^{-1}g$ at the origin has image equal to the $x$-axis.
We write our map as a composition
\[   G\times N\times G \lra E \lra \spaJ^{1} \]
where the first map is given by $(\varphi,f,\psi)\mapsto(\varphi,\varphi f\psi)$ and
the second map is given by $(\varphi,g)\mapsto g$. Clearly the second of these
maps is a fiber bundle. To understand the first map, we
fix some $(\varphi,g)\in E$. The portion of $G\times N\times G$ mapping to that
is identified with the set of all $\psi\in G$ such that $\varphi^{-1}g\psi^{-1}\in N$.
This condition on $\psi$ can also be described as saying that
the following commutes up to terms of order $\ge 5$:
\[
\CD
\RR^2 @>\psi^{-1}>> \RR^2 \\
@VV p V  @VV{p\varphi^{-1}g} V \\
\RR @>=>> \RR
\endCD
\]
where $p(x,y)=x$. If we select one such $\psi$, and we can, then all others can be obtained from
the selected one by multiplying on the left with an element of
\[  H=\{\gamma\in G~|~p\gamma=p\}, \]
a subgroup of $G$.
Hence our map $G\times N\times G \lra E$ is a
principal bundle with structure group $H$. \qed

\begin{lem}
\label{lem-slicing}
Suppose that a Lie group $L$ acts smoothly on a smooth connected manifold $M$. Let $N\subset M$
be a smooth submanifold, closed as a subset of $M$. Suppose that the restricted action map
$L\times N\to M$ is a smooth surjective submersion. Then the partition of $M$ into
$L$-orbits is locally diffeomorphic to the induced partition of $N$, multiplied
with $\RR^k$ where $k=\dim(M)-\dim(N)$.
\end{lem}

\proof Given $z\in M$, choose $(g,x)\in L\times N$ such that $gx=z$.
By assumption the differential of the action map $\alpha\co L\times N\to M$
at $(g,x)$ is a (linear) surjection $d\alpha_{(g,x)}\co T_gL\times T_xN\to T_zM$.
Its restriction to $T_xN$ is injective since it is the differential of an
embedding $N\to M$. Hence there exists a $k$-dimensional subspace $V\subset T_gL$
such that $d\alpha_{(g,x)}$ restricts to a linear isomorphism
$V\times T_xN\to T_zM$. Now choose a smooth embedding germ $s:(V,0)\to (L,g)$ such that
the differential of $s$ at $0$ is the inclusion $V\to T_gL$. Then the germ of
\[  h\co V\times N\to M~;~h(v,z)= s(v)z \]
at $(0,x)$ is a diffeomorphism germ. Clearly, for $(v_1,z_1)$ and $(v_2,z_2)$ in $V\times N$,
the elements $h(v_1,z_1)$ and $h(v_2,z_2)$ are in the same $L$-orbit if and only if
$z_1$ and $z_2$ are in the same $L$-orbit. \qed

\medskip
Putting the two lemmas together, we see that in order to understand (some of) the decomposition of $\spaJ^{1}$
into $G\times G\op$-orbits, it is sufficient to understand (some of) the induced decomposition of
$ N\subset\spaJ^{1}$. But this is already obvious from the list above. The decomposition
of the affine space $N$ can be described in terms of several linear forms
on $N$ (and a quadratic form). The linear forms are given by the coefficients $b_2,b_3,b_4,c,d_1,d_2$
of $y^2,y^3,y^4,xy,xy^2,x^2y$~, respectively. The quadratic form is $q=3b_3d_2-d_1^2$.
Let $B_2,B_3,B_4,C,Q\subset  N$
be the zero sets of $b_2,b_3,b_4,c,q$ respectively.
Now we can describe the ``relevant'' strata (intersected with $N$) as follows:
\begin{itemize}
\item \emph{fold}: $ N\smin B_2$
\item \emph{cusp}: $B_2\smin (B_3\cup C)$
\item \emph{swallowtail}: $(B_2\cap B_3)\smin(B_4\cup C)$
\item \emph{lips}: $(B_2\cap C)\smin Q$ (and $q>0$)
\item \emph{beak-to-beak}: $(B_2\cap C)\smin Q$ (and $q<0$).
\end{itemize}
The points of $ N$ which are not in any of these strata form a closed codimension 3 algebraic subset:
\[  N\smin W^{\spaJ}= (B_2\cap B_3\cap B_4)\cup(B_2\cap B_3\cap C)\cup (B_2\cap C\cap Q)~. \]

\begin{prop}
\label{prop-strat1} The complement of $ W^{\spaJ}$ in $\spaJ$ is closed, algebraic and of codimension $\ge 4$.
The stratification of $ W^{\spaJ}$ by the six strata (alias $ G\times G\op$-orbits) is in fact a filtration by
smooth submanifolds (of codimensions 0,1,2,3) as indicated in the following diagram:
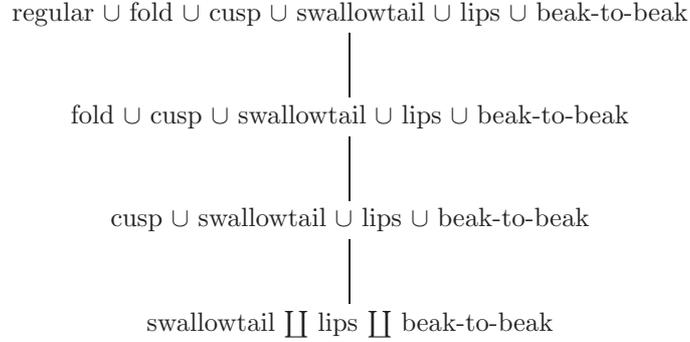
\[
\xymatrix@C=5pt{
\textup{regular $\cup$ fold $\cup$ cusp $\cup$ swallowtail $\cup$ lips $\cup$ beak-to-beak} \ar@{-}[d]  \\
\textup{fold $\cup$ cusp $\cup$ swallowtail $\cup$ lips $\cup$ beak-to-beak}  \ar@{-}[d]  \\
\textup{cusp $\cup$ swallowtail $\cup$ lips $\cup$ beak-to-beak} \ar@{-}[d] \\
\textup{swallowtail $\coprod$ lips $\coprod$ beak-to-beak}
}
\]
\end{prop}

\proof The set $\spaJ^{1}\cap  W^{\spaJ}= G( N\cap W^{\spaJ}) G$ is open in $\spaJ^{1}$ because $ N\cap W^{\spaJ}$ is open in $ N$.
Hence $ W^{\spaJ}=\spaJ^{2}\cup(\spaJ^{1}\cap  W^{\spaJ})$ is open in $\spaJ^{2}\cup\spaJ^{1}$ which in turn is open in $\spaJ$.
The same argument shows that $ W^{\spaJ}$ is algebraic in $\spaJ$, given that the two actions of $ G$ on $\spaJ$
are algebraic. The codimension of $ G( N\smin W^{\spaJ}) G=\spaJ^{1}\smin W^{\spaJ}$ in $\spaJ^{1}$ is $\ge3$
by lemma~\ref{lem-slicing}. Hence the codimension of $\spaJ^{1}\smin W^{\spaJ}$ in $\spaJ$ is $\ge4$.
The codimension of $\spaJ^{0}$ in $\spaJ$ is also $4$. \newline
The second statement follows from our analysis of the stratification of $ N$, together
with lemma~\ref{lem-slicing}. \qed

\begin{rem} {\rm All elements of $ W^{\spaJ}$ are represented
by proper maps $\RR^2\to\RR^2$ taking the origin to itself, and have a well-defined degree.
The degree is $0$ in the case of a fold
or swallowtail, but $\pm1$ in the case of a regular germ, cusp, lips or beak-to-beak.
This shows that at least four of the six strata in our stratification of $ W^{\spaJ}$ are not connected.}
\end{rem}

\section{Germs of maps from higher dimensional space to the plane} \label{sec-highleftright}
We generalize the results above by investigating
(certain) smooth map germs
\[ f\co(\RR^{n+2},0)\to (\RR^2,0) \]
for fixed $n\ge 0$.
It turns out that there is an easy reduction to the case $n=0$.

\smallskip
\subsection{Classification}
We begin with the classification up to left-right equivalence.
Again we exclude the cases where $df(0)$ has rank $0$ and note that the rank 2 case is easy.
This leaves the rank 1 case. Using appropriate linear transformations of source and target,
we may assume that
\[ df(0)=\begin{bmatrix} 0 & 0 & 0 &  \cdots & 0 & 1 & 0 \\ 0 & 0 & 0 & \cdots & 0 & 0 & 0\end{bmatrix} \]
so that the image of $df(0)$ is the $x$-axis. Writing $p\co \RR^2\to
\RR$ for the linear projection $(x,y)\mapsto x$, we can use $pf$ as
one of $n+2$ coordinates on the source and so obtain
\[  f(z_1,\dots,z_n,x,y)=(x,f_2(z_1,\dots,z_n,x,y)) \]
where $f_2\co (\RR^{n+2},0)\to (\RR,0)$ has vanishing derivative at
$0$. Then we require that the Hessian of $f_2$, restricted to
$\ker(df(0))$, be not too singular: its nullspace must have
dimension $\le 1$. (The cases where the nullspace has dimension
$\ge2$ are considered too rare to be of interest here.) There are
two cases to distinguish.

\smallskip\nin
\emph{Case 1: The nullspace of that restricted Hessian
has dimension 0.} By the
Morse lemma we can assume, after a coordinate transformation in the source (involving
only the coordinates $z_1,\dots,z_n,y$), that $f_2$ restricted to $\ker(df(0))$ is a
quadratic form alias homogeneous polynomial of degree $2$. Then $f_2$ can be viewed as a
1-parameter deformation of the restriction of $f_2$ to $\ker(df(0))$. By the classification of such deformations,
we may assume that the deformation is merely given by translations in the target
(after another coordinate transformation in the source). Then we have the form
\[  f(z_1,\dots,z_n,x,y)=(x,q(z_1,\dots,z_n,y)+g(x)) \]
where $q$ is a nondegenerate quadratic form in $n+1$ variables. Finally we may
remove the $g(x)$ term using a coordinate transformation in the target. This gives the form
\[  f(z_1,\dots,z_n,x,y)=(x,u(z_1,\dots,z_n,y)) \]
where $u$ is a nondegenerate quadratic form in $n+1$ variables. Using another
linear transformation of the source coordinates $z_1,\dots,z_n,y$ and where necessary
a reflection $(x,y)\mapsto(x,-y)$ in the target, we reduce further to the case
where $u(z_1,\dots,z_n,y)=y^2+q(z_1,\dots,z_n)$ for a quadratic form $q$ in the variables
$z_1,\dots,z_n$. Then we have the canonical form
\[  f(z_1,\dots,z_n,x,y)=(x,y^2+q(z_1,\dots,z_n)) \]
where $q$ is a nondegenerate quadratic form in the variables $z_1,\dots,z_n$.

\smallskip\nin
\emph{Case 2: The nullspace of that restricted Hessian has dimension 1.} We may assume
that the nullspace is the $y$-axis. Let
$K=\{(z_1,\dots,z_n,0,0)\}\subset \RR^{n+2}$. By the Morse lemma
applied to $f_2|_K$, we may assume that $f_2|_K$ is a nondegenerate
quadratic form (after a suitable coordinate transformation in the
source involving only $z_1,\dots,z_n$). Now we can view $f$ as a
2-parameter deformation (parameters $x$ and $y$) of $f_2|_K$. By the
classification of such deformations, we may assume that the
deformation is merely given by translations in the target. Then
\[ f(z_1,\dots,z_n,x,y)=(x,f_2^r(x,y)+q(z_1,\dots,z_n)) \]
where we write $f_2^r$ to indicate a ``reduced'' form of $f_2$. In words, $f_2$ has the
form of a function germ $f_2^r$ which only depends on the variables $x$ and $y$, and has vanishing
first derivative at $0$, plus a nondegenerate quadratic form $q$ which depends only on the
\emph{other} variables $z_1,\dots,z_n$. The second derivative at $0$ of $f_2^r$
restricted to the $y$-axis is zero, because we are not in ``case 1''.

\medskip
The analysis in case 1 above is fairly complete. We call this type a \emph{fold}. In case 2,
it is natural to proceed by imposing a condition: namely, that the
germ
\[ (x,y)\mapsto(x,f_2^r(x,y)) \]
have one of the types
\emph{cusp}, \emph{swallowtail}, \emph{lips} or \emph{beak-to-beak} described earlier in this section.
Then we get the list of normal forms
\begin{align}
\label{eqn-foldnorm} {\textit Fold}: f(z_1,\dots,z_n,x,y) & =  (x,y^2+q(z_1,\dots,z_n)) \\
\label{eqn-cuspnorm} {\textit Cusp}: f(z_1,\dots,z_n,x,y)& =  (x,y^3+xy+q(z_1,\dots,z_n))  \\
\label{eqn-swallownorm} {\textit Swallowtail}: f(z_1,\dots,z_n,x,y) & =  (x,y^4+xy+q(z_1,\dots,z_n)) \\
\label{eqn-lipsnorm} {\textit Lips}: f(z_1,\dots,z_n,x,y) & =  (x,y^3+x^2y+q(z_1,\dots,z_n)) \\
\label{eqn-beaknorm} {\textit Beaktobeak}: f(z_1,\dots,z_n,x,y) & =  (x,y^3-x^2y+q(z_1,\dots,z_n)).
\end{align}
In these formulae, $q$ is a nondegenerate quadratic form. It is easy to see that
the five types are distinguishable in coordinate free terms. For example, in the cusp and swallowtail cases,
the singularity set in the source is a smooth submanifold of dimension 1,
but in the lips and beak-to-beak cases, it is not. The cusp case can be distinguished
from the swallowtail case because the singularity sets in the target are not equivalent.

\medskip
The above reduction procedure extends easily to 1-parameter
families. Indeed, suppose that we have a smooth function germ
$(\RR\times \RR^{n+2},0)\to (\RR^2,0)$ which we want to regard as a
1-parameter family of germs
\[  f_t:(\RR^{n+2},0)\to (\RR^2,0) \]
with $t\in \RR$ in a neighborhood of $0$.
Suppose that the first derivative of each $f_t$ at $0$ has rank 1,
and also that $f_0$ has the ``reduced'' form
\[  f_0(z_1,\dots,z_n,x,y)=(x,f_{0,2}^r(x,y)+q_0(z_1,\dots,z_n)) \]
where $q_0$ is a nondegenerate quadratic form in $n$ variables.
Then there exist diffeomorphism germs
\[ \psi_t:(\RR^{n+2},0)\to(\RR^{n+2},0)~,\qquad \varphi_t:(\RR^2,0)\to (\RR^2,0) \]
depending smoothly on $t$, with $\psi_0=\id$ and $\varphi_0=\id$, such that
$\bar f_t=\varphi_tf_t\psi_t$ is in reduced form,
\[ \bar f_t(z_1,\dots,z_n,x,y)\mapsto (x,\bar f_{t,2}^r(x,y)+q_t(z_1,\dots,z_n)). \]
Here $q_t$ is a nondegenerate quadratic form in $n$ variables.
Therefore we have proved lemma~\ref{lem-highfiberbundle} below. \newline

\subsection{Codimension and stratification} \label{subsubsec-codimstrat}
Let $\spaJ$ be the finite dimensional real vector space of polynomial maps
$\RR^{n+2}\to \RR^2$ of degree $\le 4$, with vanishing constant term. We write
\[ \spaJ=\spaJ^{2}\cup \spaJ^{1}\cup\spaJ^{0} \]
where $\spaJ^{i}$ consists of the polynomials whose linear term has
rank $i$. Let $G$ be the set of polynomial maps of degree $\le 4$ from $\RR^{n+2}$ to $\RR^{n+2}$, with vanishing constant term and
invertible linear term. Under composition and truncation, $G$
becomes a group, and this group acts on the right of $\spaJ$ by
composition. Let $W^{\spaJ}\subset\spaJ$ be the union of the six
strata \emph{regular, fold, cusp, swallowtail, lips} and
\emph{beak-to-beak}. Let $D\subset \spaJ^{1}$ be the closed
subset consisting of the elements whose second ``Porteous''
derivative has a nullspace of dimension $>1$. Let $\sF$ be the
space of nondegenerate quadratic forms in $n$ real variables
$z_1,\dots,z_n$. We write $G_{ol}$ for the old $ G$ of lemma~\ref{lem-fiberbundle}, and
$ N_{ol}$ for the old $ N$ of lemma~\ref{lem-fiberbundle}.

\begin{lem}
\label{lem-highfiberbundle}
The restricted action map
\begin{eqnarray*}
\,\, G_{ol}\times N_{ol}\times\sF\times  G & \lra & \spaJ^{1}\smin D \\
(\varphi,f,q,\psi) & \mapsto & \varphi(f+q)\psi~,
\end{eqnarray*}
where $f+q$ is shorthand for the map
\[  (z_1,\dots,z_n,x,y)\mapsto (x,f_2(x,y)+q(z_1,\dots,z_n))~, \]
is a surjective submersion. \qed
\end{lem}

\medskip
This puts us in a position to use lemma~\ref{lem-slicing}. Hence the
partition of $\spaJ^{1}\smin D$ into $G_{ol}\times G\op$ orbits is locally
diffeomorphic to the induced partition of $ N_{ol}\times\sF$. But the latter is
essentially the partition of $ N_{ol}$ into $ G_{ol}\times G_{ol}\op$ orbits
multiplied with a certain partition of $\sF$ where each part is a union of path components.

\begin{cor}
\label{cor-highstrat1} The complement of $W^{\spaJ}$ in $\spaJ$ is closed, algebraic
and of codimension $\ge n+4$.
The stratification of $ W^{\spaJ}$ by the six strata is given by a nested sequence of
smooth algebraic subvarieties of $W^{\spaJ}$ of codimensions $0$, $n+1$, $n+2$, $n+3$, respectively, as
indicated in the following diagram:
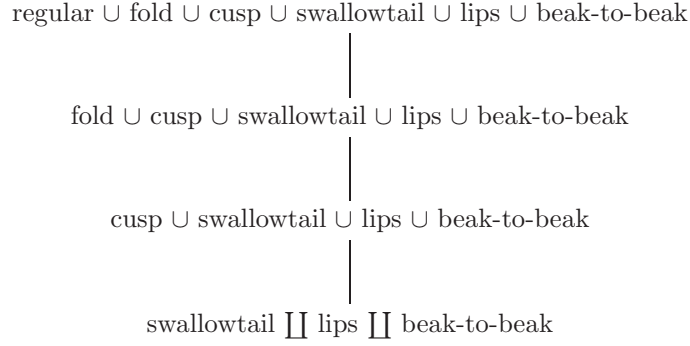
\[
\xymatrix@C=5pt{
\textup{regular $\cup$ fold $\cup$ cusp $\cup$ swallowtail $\cup$ lips $\cup$ beak-to-beak} \ar@{-}[d]  \\
\textup{fold $\cup$ cusp $\cup$ swallowtail $\cup$ lips $\cup$ beak-to-beak}  \ar@{-}[d]  \\
\textup{cusp $\cup$ swallowtail $\cup$ lips $\cup$ beak-to-beak} \ar@{-}[d] \\
\textup{swallowtail $\coprod$ lips $\coprod$ beak-to-beak}
}
\]
It is invariant under the action of $\,\, G_{ol}\times G\op$.
The ``regular'' stratum is a single orbit of that action. The ``fold'' stratum
falls into $\lfloor n/2+3/2\rfloor$ orbits, and
the ``cusp'', ``swallowtail'', ``lips'' and ``beak-to-beak''
strata fall into $\lfloor n/2+1\rfloor$ orbits each.
\end{cor}

\proof Most of this has already been established. The left-right equivalence class counts are
obtained by counting components of suitable spaces
of nondegenerate quadratic forms, modulo sign change. In the fold case, we have to look at
nondegenerate quadratic forms in $n+1$ variables. The components are classified by the signature,
which can be $n+1$, $n-1$, \dots, $-n-1$.
If we allow sign change, as we must, only the absolute value
of the signature remains, so there are $\lfloor n/2+3/2\rfloor$ types. In the remaining cases, we are
looking at nondegenerate quadratic forms in $n$ variables. There are $\lfloor n/2+1\rfloor$ types.
\qed

\section{Multigerms of maps}
Let $L$ be a smooth manifold and $S\subset L$ a finite nonempty subset. We are interested in \emph{multigerms} of smooth
maps $f\co (L,S)\to (\RR^m,0)$. Such a multigerm is, strictly speaking, an equivalence class of pairs $(U,f)$ where $U$ is a neighborhood of $S$ in $L$
and $f\co U\to \RR^m$ is a smooth map taking all of $S$ to $0$. Two such pairs $(U_0,f_0)$ and $(U_1,f_1)$ are
equivalent if $f_0$ and $f_1$ agree on some neighborhood of $S$ contained in $U_0\cap U_1$. \newline
The germs $(L,s)\to (\RR^m,0)$ for $s\in S$, obtained by restriction or localization from $f$, are the \emph{branches} of
the multigerm $f\co (L,S)\to (\RR^m,0)$. Consequently $|S|$ is the \emph{number of branches}.

\begin{defn} \label{defn-leftrightmultigerm}
{\rm Two multigerms $f\co(L,S)\to (\RR^m,0)$ and $g\co (L',S')\to (\RR^m,0)$ are \emph{left-right equivalent} if
there exist a diffeomorphism germ $\psi\co(L,S)\to (L',S')$, extending some bijection $S\to S'$, and a diffeomorphism germ
$\sigma\co(\RR^m,0)\to(\RR^m,0)$ such that $g=\sigma f\psi^{-1}$. The multigerms $f$ and $g$ are \emph{right equivalent}
if there exists a diffeomorphism germ $\psi\co(L,S)\to (L',S')$, extending some bijection $S\to S'$, such that $g=f\psi^{-1}$.
}
\end{defn}

There are similar notions of left-right equivalence and right equivalence
for \emph{multijets}. We have in mind the finite set $S_r=\{1,2,\dots,r\}$
for $r\ge 1$, and two elements $f,g$ of
\begin{equation} \label{eqn-multijet} \prod_{x\in S_r} \spaJ \end{equation}
where $\spaJ$ is the vector space of polynomial mappings of degree $\le z$ from $\RR^\ell$ to $\RR^m$, with vanishing
constant term, for some $z>0$. (Soon we will take $z=4$ or $z\ge 4$ and $\ell=n+2$, $m=2$ as in previous sections.)
If necessary, we refer to $z$ as the \emph{order} of the multijet while $r$ is (still) the \emph{number of branches}.

\begin{defn} \label{defn-leftrightmultijet}
{\rm  The multijets $f$ and $g$ (of order $z$) are \emph{left-right equivalent} if there
exist jets (of order $z$) of diffeomorphisms $\psi$ from $(S_r\times\RR^\ell,S_r)$ to $(S_r\times\RR^\ell,S_r)$,
extending some permutation of $S_r$~, and
$\sigma$ from $(\RR^m,0)$ to $(\RR^m,0)$, such that $g=\sigma f\psi^{-1}$.
(We have identified $S_r$ with $S_r\times\{0\}\subset S_r\times\RR^{n+2}$.) The multijets $f$ and $g$ are
\emph{right equivalent}
if there is a jet (of order $z$) of diffeomorphisms $\psi$ from $(S_r\times\RR^\ell,S_r)$ to $(S_r\times\RR^\ell,S_r)$,
extending some permutation of $S_r$~, such that $g=f\psi^{-1}$.
}
\end{defn}

\begin{rem} \label{rem-sectionplan} {\rm
For the rest of the section we take $\RR^m=\RR^2$ as the target manifold and focus on source manifolds $L$
of dimension $n+2$, unless otherwise stated. Our goal is to select for each $r\ge 1$
an open semi-algebraic subset $\mathfrak{X}_r\subset \prod_{x\in S_r} \spaJ$~, closed under
\emph{right equivalence}, in such a way that a number of desirable conditions are satisfied.
The multijets which belong to $\mathfrak{X}_r$~, for some $r$, and the multigerms $(L,S)\to(\RR^2,0)$ whose multijets
belong to $\mathfrak{X}_r$ (in multilocal coordinates about $S\subset L$, assuming that $S$ admits a bijection to $S_r$)
are called \emph{admissible}. Among the desirable conditions is
\begin{itemize}
\item[(a)] Naturality: for  an admissible multigerm
from $(L,S)$ to $(\RR^2,0)$, and any nonempty subset $T$ of $S$,
the induced multigerm from $(L,T)$ to $(\RR^2,0)$ is admissible.
More generally, for any admissible multigerm $f$ from $(L,S)$ to $(\RR^2,0)$,
there exists a neighborhood $U$ of $S$ in $L$ with the following property.
For any finite nonempty subset $T$ of $U$ such that $f|T$ is constant, the multigerm of $f$ at $T$,
minus that constant, is admissible.
\end{itemize}
Suppose that (a) holds and let $f\co L\to \RR^2$ be a smooth map, where $\dim(L)=n+2$. We say that
$f$ is admissible if, for every finite nonempty subset $S\subset L$ such that $f|S$ is constant,
the multigerm of $f$ at $S$, minus that constant, is admissible.
Conditions (b) and (c) below ensure, loosely speaking, that for $L$ as
above the cohomology of the space of admissible smooth maps $L\to \RR^2$ admits a description
in terms of the cohomology of the spaces of admissible smooth multigerms $(L,S)\to(\RR^2,0)$,
where $S$ runs through the finite nonempty subsets of $L$.
(We will not explain here \emph{how} conditions (b) and (c) ensure that; see \cite{RWIII} instead.)
For finite nonempty $S\subset L$ and a non-admissible germ
\[ g\co (L,S) \lra (\RR^2,0), \]
a nonempty subset $T$ of $S$ is \emph{a minimal bad event} if the multigerm
$(L,T)\to(\RR^2,0)$ obtained from $g$ by restriction
is non-admissible and $T$ has no proper nonempty subset with the same property. A nonempty
subset $T$ of $S$ is a \emph{bad event} for $g$ if it is a union of minimal bad events for $g$.
The \emph{size} of $g$ is the maximum cardinality of a bad event for $g$. The \emph{complexity}
of $g$ is the maximum of the integers $k$ such that there exists a chain
of bad events
$T_0\subset T_1\subset \dots \subset T_{k-1}\subset T_k$
where $T_i\ne T_{i+1}$ for $i=0,\dots,k-1$.
\begin{itemize}
\item[(b)] The codimension $c_*(s)$ of the set of non-admissible multigerms of size $s$ is at least $sn+4$.
\item[(c)] For $k\le s$, the codimension $c_*(s,k)$ of the set of non-admissible multigerms of size $s$ and complexity $k$ satisfies
\[  \lim_{s\to\infty} (c_*(s,k)-sn-k)~=~\infty~. \]
\end{itemize}
More precisely: in the multijet space ~(\ref{eqn-multijet}), the subset of non-admissible
multijets of size $s$ and complexity $k$ (where $k\le s\le r$) is a semi-algebraic subset, with a minimum codimension
which we denote by $c_*(s,k,r)$. Let $c_*(s,k)= \min_r \{c_*(s,k,r)\}$.
It is easy to see that $c_*(s,k)=c_*(s,k,s)$. Let $c_*(s)=\min_k \{c_*(s,k)\}$.
These definitions of codimension should be used in conditions (b) and (c).
See also remark~\ref{rem-othercodim}.
}
\end{rem}

\begin{rem} \label{rem-othercodim} {\rm
Let $X$ be the vector space of all smooth maps to $\RR^2$ from a
smooth $(n+2)$-manifold $L$, closed for simplicity. In $X\times L\times\dots\times L$, form
the subset of all $(f,x_1,\dots,x_s)$ such that $x_1,\dots,x_s$ are distinct while
$f(x_1)=\cdots= f(x_s)=:a~$, and $S=\{x_1,\dots,x_s\}$ is a bad event of complexity $k$ (and size $s$) for the
multijet of $f-a$ at $S$. Multijet transversality theorems imply that this subset has a well defined
minimum codimension which turns out to be
\[  c(s,k):=c_*(s,k)+2(s-1)~. \]
It is therefore tempting to think, but not obviously meaningful,
that the subset of $X$ consisting of all non-admissible $f$ which have
some bad event of size $s$ and complexity $k$ has codimension at least
\[  C(s,k):=c(s,k)-s(n+2)=c_*(s,k)-sn-2~. \]
We justify this idea in \cite{RWIII}, following Vassiliev.
Now condition (c) of remark~\ref{rem-sectionplan} implies
\[ \lim_{s\to\infty} (C(s,k)-k)~=~\infty \]
and the inequality in condition (b) implies $C(s)\ge 2$. These are the properties that we are after.
}
\end{rem}

\bigskip
We now describe our subsets $\mathfrak{X}_r\subset \prod_{x\in S_r} \spaJ$\,, taking $z\ge 4$. Later we
point out that the $\mathfrak{X}_r$ for all $r\ge 0$ together constitute a minimal choice, under the conditions
listed in remark~\ref{rem-sectionplan} and an additional condition described in lemma~\ref{lem-2span}
and remark~\ref{rem-2span}. 
In an earlier version of this article, the additional condition was that $\mathfrak X_r$ should be closed
under left-right equivalence for each $r$. This turned out to be insufficient for a characterization
of the sets $\mathfrak X_r$ by minimality. 

\begin{defn} \label{defn-admissible} {\rm A multijet
\[ (f_x)_{x\in S_r} \in \prod_{x\in S_r} \spaJ \]
is admissible, i.e., is an element of $\mathfrak X_r$~, if and only if
\begin{itemize}
\item[--] each $f_x$
belongs to one of the types \emph{regular, fold, cusp, swallowtail, lips, beak-to-beak};
\item[--] at most one of the $f_x$ does not belong to one of the types \emph{regular, fold};
\item[--] \emph{either} for all singular $f_x$~, the images of their linear parts are distinct elements of $\RR P^1$~;
\item[--] \emph{or} all singular $f_x$ are of type \textit{fold}, and for precisely two of them the images of their
linear parts agree; in that case the two fold curves in the target make an ordinary (first order) tangency at the origin.
\end{itemize}
}
\end{defn}

\medskip 
From the definition, $\mathfrak X_r$ decomposes into manifold strata with names such as \emph{one cusp and $r-1$ folds, 
making $r$ distinct directions in target}, or \emph{two kissing folds and $(r-2)$ other folds, making $r-1$ 
distinct directions in the target}.

\begin{expl}\label{expl-minbad} {\rm Let $f=(f_x)_{x\in S_r}$ be a multijet
and let $T\subset S_r$ be a minimal bad event
for $f$. If $T=\{x\}$ has cardinality 1 then
\begin{itemize}
\item[(i)] $f_x$ is a jet
which is not of type {\it fold, cusp, swallowtail, lips} or {\it beak-to-beak}.
\end{itemize}
If $T=\{x,y\}$ is of cardinality 2,
then $f_x$ and $f_y$ are both of type {\it fold, cusp, swallowtail, lips} or {\it beak-to-beak},
and one of the following applies:
\begin{itemize}
\item[(ii)] neither $f_x$ nor $f_y$ are of {\it fold} type;
\item[(iii)] exactly one of the two is of {\it fold} type and the image of the linear part is the same for both;
\item[(iv)] both are of fold type and their fold lines make a higher tangency (double, triple etc.) in the target.
\end{itemize}
If $T=\{x,y,z\}$ has cardinality 3, then $f_x$, $f_y$ and $f_z$ are all of type {\it fold,
cusp, swallowtail, lips} or {\it beak-to-beak}, and one of the following applies:
\begin{itemize}
\item[(v)] exactly one of $f_x$, $f_y$, $f_z$ is not of fold type, with image of
differential $\ell$, while the other two are folds and
share the image $\ell'$ of their linear part, making an ordinary tangency in the target, $\ell'\ne \ell$;
\item[(vi)] $f_x$, $f_y$ and $f_z$ are all of {\it fold} type, the image of the linear part is the
same for all, and any two make an ordinary tangency in the target.
\end{itemize}
This covers all cases. So a minimal bad event has cardinality at most $3$. Each of the above six cases defines
a semi-algebraic subset of multijet space~(\ref{eqn-multijet}), for the appropriate $r\in \{1,2,3\}$. It is
easy to show that the
codimension is bounded below by $n+4$ in case (i), by $2n+4$ in cases (ii),(iii) and (iv), and by $3n+5$ in cases
(v) and (vi).
}
\end{expl}

\begin{defn}\label{defn-badchains} {\rm
Let $T_\bullet := T_0 \subset T_1\subset\cdots\subset T_k$
be a chain of finite nonempty sets and proper inclusions.
We define
\[ Y(T_\bullet)\subset \prod_{s\in T_k}\spaJ \]
to consist of all elements $h$ such that $T_j$ is a bad event
for $h$, for $0\le j\le k$, and there is no bad event for $h$ strictly between $T_j$ and $T_{j+1}$~, for $0\le j\le k-1$.
}
\end{defn}

\begin{lem}\label{lem-badchainscodim}
The codimension of the semialgebraic set $Y(T_\bullet)$ in $\prod_{s\in T_k}\spaJ$ is at least $|T_k|n+2k+4$ everywhere.
\end{lem}

\proof We proceed by induction on $k$.
The case where $k=0$ has been dealt with in example~\ref{expl-minbad}. In the case $k>0$,
let $T'_{\bullet}$ be the truncated chain
\[ T_0 \subset T_1\subset\cdots\subset T_{k-1}~. \]
Let $R=T_k\smin T_{k-1}$. There is a projection
\begin{equation} \label{eqn-allrestrict}
\prod_{s\in T_k} \spaJ \,\,\lra \prod_{s\in T_{k-1}} \spaJ
\end{equation}
which induces a projection
\begin{equation} \label{eqn-badrestrict}
Y(T_\bullet)\to Y(T'_\bullet).
\end{equation}
Fix some $h\in Y(T'_\bullet)$.
The fiber $F_h$ of~(\ref{eqn-badrestrict}) over $h$ is a semialgebraic subset of the fiber $E_h$ 
of~(\ref{eqn-allrestrict}) over $h$, where $E_h$ is a vector space,
\[  E_h \cong \prod_{s\in R} \spaJ ~. \]
Now it is enough to show that the codimension of $F_h$ in $E_h$ is at least $|R|n+2$. In the case where
$|R|>1$ this is instantly clear. Indeed the codimension of $F_h$ in $E_h$ is at least $|R|(n+1)$, because
the germs $g_x$ for $g\in F_h$ and $x\in R$ are all singular, and the singular subset of $\spaJ$ has codimension $n+1$.
Suppose then that $R$ is a singleton, $R=\{x\}$. Write $F_h$ as a union of three semialgebraic subsets,
one containing the elements $g$ for which $R$ is a minimal bad event, the second one containing the
elements $g$ for which $R$ participates in a minimal bad event of cardinality 2, and the last one containing
the elements $g$ for which $R$ participates in a minimal bad event of cardinality 3. Looking at the three cases
separately, we see that the jet $g_x$ for $g\in F_h$ is either singular and not of {\it fold} type, or it is of
{\it fold} type but the direction of the fold line in the target is prescribed by $h$ up to finite choice.
Hence the codimension of $F_h$ in $E_h\cong \spaJ$ is at least $n+2=|R|n+2$. \qed

\begin{thm} The subsets $\mathfrak{X}_r$ of definition~{\rm \ref{defn-admissible}}
together satisfy conditions {\rm (a), (b)} and {\rm (c)}
of remark {\rm \ref{rem-sectionplan}}.
\end{thm}

\proof By inspection, condition (a) is satisfied. For condition (b), let
\[    Z_r \subset \prod_{x\in S_r} \spaJ \]
consist of all the multijets $f=(f_x)$ such that all of $S_r$ is a bad event for $f$. We need to show that
the codimension of $Z_r$ in $\prod_x\spaJ$ is at least $rn+4$. Let
\[  Q_r\subset \prod_{x\in S_r} \spaJ \]
consist of all the $f=(f_x)$ such that $f_x$ is singular for every $x\in S_r$. Then $Z_r\subset Q_r$
and the codimension of $Q_r$ in $\prod_x \spaJ$ is $r(n+1)$. Therefore it is enough to show that the codimension of $Z_r$
in $Q_r$ is at least 1 when $r=3$, at least $2$ when $r=2$ and at least $3$ when $r=1$. That is
easily done by inspection. \newline
Next we verify condition (c). We look for lower bounds for $c_*(s,k,r)$ since $c_*(s,k)= \min_r \{c_*(s,k,r)\}$.
It is understood that $k\le s\le r$.
If there are no non-admissible multijets of size $s$ and complexity $k$ in $\prod_{x\in S_r}\spaJ$~, then
\[  c_*(s,k,r)= r\cdot\dim(\spaJ)> s(n+2)=sn+2s~. \]
If there are such multijets, then $k+1\ge s/3$ because minimal bad events have cardinality $\le 3$. By
lemma~\ref{lem-badchainscodim}, we have
\[  c_*(s,k,r)\ge sn+2k+4  \]
so that $c_*(s,k,r)-sn-k\ge k+4> s/3$. Therefore 
\[ c_*(s,k)-sn-k > s/3 \]
which establishes condition (c). \qed

\bigskip
Suppose that $f^{(t)}\co (L,S) \to (\RR^m,0)$ are multigerms, depending smoothly on $t\in[0,1]$. Here the dimension 
of $L$ is arbitrary. We say that the family
$(f^{(t)})$ is \emph{left-right trivial} if there exist diffeomeorphism germs
$\psi^{(t)}\co (L,S)\to (L,S)$ and $\sigma^{(t)}\co (\RR^m,0)\to(\RR^m,0)$, depending smoothly on $t\in[0,1]$,
such that $\psi^{(0)}=\id$, $\sigma^{(0)}=\id$ and 
\[ f^{(t)}\psi^{(t)}=\sigma^{(t)}f^{(0)}. \]

\begin{defn} \label{defn-qspan} {\rm Let $q$ be a positive integer. Two
multigerms $f\co (L,S)\to (\RR^m,0)$ and $g\co (L',S') \to (\RR^m,0)$ are \emph{q-span left-right equivalent} if there
exists a family of multigerms $\big(f^{(t)}\co (L,S)\to (\RR^m,0)\big)$, 
depending smoothly on $t\in[0,1]$, such that
\begin{itemize}
\item $f^{(0)}=f$ and $f^{(1)}$ is left-right equivalent to $g$
\item for every nonempty $T\subset S$ of cardinality $\le q$, the family of multigerms
$\big(f^{(t)}\co (L,T)\to (\RR^m,0)\big)$ is left-right trivial.
\end{itemize}
}
\end{defn}

There is a similar definition for multijets. Note that left-right equivalence (for multigerms or multijets)
implies $q$-span left right equivalence, and the two notions are identical for multigerms or multijets
with branch number $\le q$. \newline
Our reasons for making such a definition is that left-right equivalence for multigerms and multijets with large
branch number $r$ is hard to handle. By contrast, $q$-span left-right equivalence for multigerms and multijets with 
branch number $r$ is as manageable as left-right equivalence for multigerms 
with branch number not greater than $q$.
See \cite{RWV} for calculations and illustrations.
We mention just one simple but striking example. For fixed $\ell=\dim(L)\ge 2$,
there are uncountably many left-right equivalence classes of multigerms $f\co(L,S)\to (\RR^2,0)$
such that $|S|=4$, all branches of $f$ are fold singularities, and the fold curves make four distinct
directions in the target (at the origin). But if these multigerms are classified
by $2$-span or even $3$-span left-right equivalence, then there are only
finitely many equivalence classes.

\begin{lem} \label{lem-2span} The sets $\mathfrak X_r$
are closed under 2-span left-right equivalence.
\end{lem}

\proof The key observation is that each stratum of $\mathfrak X_r$ for $r>2$ can be 
characterized in terms of preimages of strata in $\mathfrak X_2$ under various 
projections, while $\mathfrak X_2$ is obviously closed under left-right equivalence. For example, a 
multijet in $\prod_{x\in S_5}P_*$ is of  
type \emph{one swallowtail and four folds, making five distinct directions in the target} if and only if it has 
the following census of sub-multijets with two branches: four of type \emph{one swallowtail and one fold, making 
distinct directions in the target} and six of type \emph{two folds making distinct directions in the target}. \qed

\begin{rem} \label{rem-2span} {\rm 
(i) Each stratum of $\mathfrak X_r$ is a union of finitely many $2$-span left-right equivalence classes
which are open and closed in the stratum. The equivalence classes making up each stratum can be distinguished 
by quadratic form data, roughly as in corollary~\ref{cor-highstrat1}. We omit the details and refer 
to \cite{RWV} for the necessary calculations, which are unexciting in any case. \newline 
(ii) The sets $\mathfrak X_r$ are minimal if we insist on conditions (a), (b) and (c)
in remark~\ref{rem-sectionplan} and the property of being closed under $2$-span left-right equivalence.
Instead of giving a proof, we give a few examples to explain
why the sets $\mathfrak X_r$ have to be as big as they are. We are dealing with
multigerms $f\co (L,S)\to (\RR^2,0)$ where $\dim(L)=n+2$. \newline
Suppose to start with that $S=\{1,2\}$, that the first branch of $f$ is a cusp,  
the second is a fold, and the two branches determine distinct directions (elements of $\RR P^1$) in
the target. The left-right equivalence class of $f$ is a subset $Y_2$ of the
multijet space. As such it has codimension $2n+3$, of which $n+2$
are contributed by the cusp and $n+1$ by the fold. 
Therefore by our conditions on $\mathfrak X_2$~, specifically condition (b) in~\ref{rem-sectionplan}, we must have
$Y_2\subset \mathfrak X_2$. (A similar but easier argument shows that the multijets with branch number 2 made up of
two fold singularities, distinct directions in the target, are all in $\mathfrak X_2$~.) \newline
Suppose next that $S=\{1,2,3\}$ where $r\ge 2$, that the first branch of $f$ is a cusp,
the other two branches are folds,
and the three branches determine three distinct directions (elements of $\RR P^1$) in the target. The 2-span
left-right equivalence class $Y_3$ of $f$
has codimension $3n+4$ in the multijet space ($n+2$
contributed by the cusp and $n+1$ by each fold). 
We do not violate condition (b) by declaring that $Y_3$ is in the complement of $\mathfrak X_3$.
Let us try to make such a declaration and see whether we run into problems. \newline
In order to see some problems, we look at multijets in $\prod_{x\in S} P_*$ where $S=S_r=\{1,2,3,\dots,r\}$,
with $r>3$, where the first branch is a cusp and the other branches are folds, all making
distinct directions in the target. The 2-span
left-right equivalence class $Y_r$ of $g$
has codimension $rn+r+1$ in the multijet space. It is easy to construct $g$
in such a way that for each subset $T$ of $S$ of the form
$T=\{1,2,t\}$ with $3\le t\le r$, the multijet obtained by (multi-)localization at $T$
is in $Y_3$~, therefore not in $\mathfrak X_3$. For the multijet $g$ itself, 
each of the subsets $\{1,2,t\}$ for $3\le t\le r$ is then a minimal bad event and so the subsets
$\{1,2,\dots,t\}$ for $3\le t\le r$
are bad events. So the complexity $k$ of $g$ is at least $r-3$ and the size $s$ is $r$. The same
must be true for all multijets in $Y_r$\,.
We calculate
\[ (rn+r+1) - sn-k \le  (rn+r+1)-rn-(r-3)= 4~. \]
This does not tend to infinity when $s=r$ tends to infinity. Therefore condition (c) in remark~\ref{rem-sectionplan}
is violated. This contradiction proves that $\mathfrak X_3$ must contain $Y_3$\,.
A similar argument by contradiction, using $Y_2\subset \mathfrak X_2$ and $Y_3\subset \mathfrak X_3$,
proves that $Y_4\subset \mathfrak X_4$\,. Similar arguments by contradiction show that all multigerms
of the form $f\co (L,S)\to (\RR^2,0)$ where one branch is a cusp, the other ones are folds, and all make
distinct directions in the target, have their multijets in $\mathfrak X_r$ where $r=|S|$.
}
\end{rem}

\section{Appendix: Basic results from singularity theory}
\label{sec-prelim}
We rely mostly on the excellent book by Martinet \cite{Martinet} for definitions and theorems.
Another very readable text is \cite{Broecker}, but that is exclusively concerned with
singularities of functions (target $\RR$).

\smallskip
We take the definition of an \emph{unfolding} of a smooth map germ $(\RR^s,0)\to(\RR^t,0)$
from \cite[ch.XIII]{Martinet}.
For the definition of an \emph{isomorphism} between two unfoldings (of the same map germ, and with the
same parameter space) we also rely on the same source. Note that \cite{Broecker} has a definition
(in the case $t=1$) which is slightly more restrictive in some respects, but less restrictive in other
respects because it allows for a change of the parameter space.

\smallskip
Following \cite{Martinet}, we call an unfolding $F$ (with parameter space $\RR^p$) of a smooth map germ $f$
\emph{universal} if every other unfolding (with parameter space $\RR^q$, say) of $f$
is isomorphic to $h^*F$ for some germ $h\co (\RR^q,0)\to (\RR^p,0)$. For a universal
$F$ with minimal parameter space dimension $p$~, Martinet uses the expression \emph{minimal universal},
which we shorten to \emph{miniversal}. (Br\"{o}cker uses instead \emph{versal} for Martinet's universal,
and \emph{universal}
for Martinet's minimal universal.)

\begin{defn} {\rm
Let $\sE_{s,t}$ be the real vector space of all smooth map germs from
$(\RR^s,0)$ to $\RR^t$.
In the case $t=1$, we write $\sE_s$ instead of $\sE_{s,t}$. In the general case,
$\sE_{s,t}$ is a module over the ring $\sE_s$ by
$(u\cdot g)(x)=u(x)\cdot g(x)$
for $u\in \sE_s$ and $g\in \sE_{s,t}$.}
\end{defn}

\smallskip
\begin{defn}
\label{defn-tgtspacemap}
{\rm
The \emph{tangent space} $Tf$
of a germ $f\co(\RR^s,0)\to(\RR^t,0)$ is the vector subspace
\begin{equation}  \{ df\cdot X+Y\circ f \}\,\,\subset\,\, \sE_{s,t}
\end{equation}
where $X$ and $Y$ run through all the vector field germs defined near the origin
on $\RR^s$ and $\RR^t$~, respectively, and $df$ is the total derivative of $f$. The tangent
space is typically not an $\sE_s$ submodule.
But it is an $\sE_t$ submodule of $\sE_{s,t}$ for the action of $\sE_t$ on $\sE_{s,t}$ defined
in terms of $f$ by
\begin{equation}  (u\cdot g)(x)=u(f(x))\cdot g(x)
\end{equation}
for $u\in\sE_t$ and $g\in \sE_{s,t}$.}
\end{defn}

\begin{thm}
\label{thm-unfold}
\emph{(Main theorem on unfoldings.)}
An unfolding
\[ F\co(\RR^p\times\RR^s,0)\to(\RR^p\times\RR^t,0) \]
of a germ $f\co(\RR^s,0)\to(\RR^t,0)$ is universal if and only if
the differential at $0$ of the adjoint
$F^\ad\co(\RR^p,0)\to \sE_{s,t}$
is transverse to $Tf$. \qed
\end{thm}

\noindent\emph{Remark.} We have not specified a norm on $\sE_{s,t}$. Nevertheless,
$F^\ad$ has a well defined differential at $0$, the linear map $dF^\ad(0)\co\RR^p\to \sE_{s,t}$
defined by
\begin{align}
v & \mapsto  \left(x\mapsto \lim_{t\to 0}\frac{F(tv,x)-F(0,x)}{t}\right)
\end{align}
for $v\in \RR^p$ and $x\in \RR^s$, with $x$ sufficiently close to $0$.
The transversality condition means that $\,\,\im(dF^\ad(0))+Tf=\sE_{s,t}$.

\begin{cor} Let $f\co (\RR^s,0)\to(\RR^t,0)$ be a germ such that
$Tf$ has finite codimension in $\sE_{s,t}$. Suppose that
\[ g^{(1)},\dots,g^{(p)}\in \sE_{s,t} \]
generate $\sE_{s,t}/Tf$ as a vector space. Then $F\co (\RR^p\times\RR^s,0)\to(\RR^p\times\RR^t,0)$
defined by
\begin{equation} (z,x)\mapsto f(x)+\sum_i z_ig^{(i)}(x) \end{equation}
is a universal unfolding of $f$. \qed
\end{cor}

\begin{lem}
\label{lem-unfoldinvariant}
Let $F,G\co(\RR^p\times\RR^s,0)\to(\RR^p\times\RR^t,0)$ be unfoldings
of a germ $f\co(\RR^s,0)\to(\RR^t,0)$. If $F$ and $G$ are isomorphic as unfoldings
of $f$, then the linear map
\[  dF^\ad(0)-dG^\ad(0)\co \RR^p\lra \sE_{s,t} \]
factors through $Tf\subset \sE_{s,t}$. \qed
\end{lem}

\medskip
\noindent\emph{Remark.} This means that the composition
\begin{equation}
\CD
\RR^p @> dF^\ad(0)>> \sE_{s,t} @>\textup{proj.}>> \sE_{s,t}/Tf
\endCD
\end{equation}
is an \emph{isomorphism invariant} of the unfolding $F$ (of a fixed germ $f$,
and with fixed parameter space $\RR^p$).

\medskip
We conclude this section with a few calculations of tangent spaces of germs, in increasing order of difficulty.
These are used in section~\ref{sec-lowleftright}.

\smallskip
\begin{expl}
\label{expl-fold} {\rm
Let $f\co(\RR^2,0)\to (\RR^2,0)$ be the germ given by
\[  f(x,y)=(x,y^2). \]
This is one of the germs shown to be stable by Whitney in his investigation
of singularities of maps from the plane to the plane. Stable germs have trivial
miniversal unfoldings; equivalently, $Tf=\sE_{2,2}$. It is also easy to verify
by direct calculation that $Tf=\sE_{2,2}$.
}
\end{expl}

\begin{expl}
\label{expl-cusp} {\rm
Let $f\co(\RR^2,0)\to (\RR^2,0)$ be the germ given by
\[  f(x,y)=(x,y^3-xy). \]
This is again one of Whitney's stable germs. Therefore
$Tf=\sE_{2,2}$ and the miniversal unfolding of $f$ is trivial.
\newline As an alternative, here is a direct proof of $Tf=\sE_{2,2}$
using the Mather-Malgrange preparation theorem. We view
$\sE_{2,2}=\sE_{s,t}$ as a module over $\sE_t=\sE_2$ as in
definition~\ref{defn-tgtspacemap}. We have
$\sM_t\sE_{s,t}=\{f_1\cdot g+f_2\cdot h~|~g,h\in\sE_{s,t}\}$, where
the multiplication dot means ordinary multiplication of
vector-valued functions by scalar functions. Therefore
$\sE_{s,t}/\sM_t\sE_{s,t}$ has vector space dimension 6, and is
spanned by the (cosets of) the six maps
\[ (x,y)\mapsto \left\{
\begin{array}{l} (1,0) \\  (0,1) \\ (y,0) \\ (0,y) \\ (y^2,0) \\ (0,y^2).
\end{array} \right.
 \]
By the preparation theorem, these six maps generate $\sE_{s,t}$ as an
$\sE_t$ module. A slightly tedious verification shows
that they are all in the $\sE_t$-submodule $Tf$ of $\sE_{s,t}$.
Therefore $Tf=\sE_{s,t}$.
}
\end{expl}

\begin{expl}
\label{expl-unfoldlips} {\rm Let
$f\co(\RR^2,0)\to (\RR^2,0)$ be the germ given by
\[  f(x,y)=(x,y^3+x^2y). \]
Let $W\subset \sE_{s,t}=\sE_{2,2}$ be the linear subspace consisting of
all $k=(k_1,k_2)$ such that the first derivative of $y\mapsto k_2(0,y)$ at $y=0$ vanishes. This is
clearly an $\sE_t$-submodule of $\sE_{s,t}$, and it contains $Tf$. We want
to show that $Tf=W$. \newline
We have the standard description
\[ Tf=Jf+\tau f=\sE_s\{\,(1,2xy), (0,3y^2+x^2)\,\}+\sE_t\{\,(1,0),(0,1)\,\}, \]
where $\sE_s\{...\}$ and $\sE_t\{...\}$ denote the $\sE_s$ and $\sE_t$ submodules,
respectively, generated by the elements listed between the brackets.
A two-fold application of \cite[XV.2.1]{Martinet} proves that
\begin{equation}\label{2.7.09:1}  Tf+\sE_t\{\,(0,y)\,\}=\sE_{2,2}
\end{equation}
where $(0,y)$ is short for the map $(x,y)\mapsto(0,y)$. In more detail,
we know from theorem~\ref{thm-unfold} that
$F(x,y,u)=(x,y^3+x^2y+uy)$
defines a universal (not miniversal) unfolding, with two unfolding parameters $x$ and $u$,
of the germ $y\mapsto y^3$. By \cite[XV.2.1]{Martinet} it follows that $F$
is a stable germ. But $F$ is also a one-parameter unfolding of the germ $f$.
Then \cite[XV.2.1]{Martinet} can be applied in the opposite
direction, which leads to equation~(\ref{2.7.09:1}). \newline
Hence it is enough to check that $\sM_t\cdot(0,y) \subset Tf$.
As $Tf$ is an $\sE_t$-module, that reduces to showing that
\begin{align*} (0,xy) & \in  Tf \\
(0,y^4+x^2y^2) & \in  Tf~.
\end{align*}
For the first of these, write $2(0,xy)=(1,2xy)-(1,0)$
where $(1,2xy)\in Jf$ and $(1,0)\in \tau f$. For the second, write
\[ 9(0,y^4+x^2y^2)=3y^2(0,3y^2+x^2)+2x^2(0,3y^2+x^2)-2x^4(0,1) \]
where $3y^2(0,3y^2+x^2)\in Jf$ and $2x^2(0,3y^2+x^2)\in Jf$ and $2x^4(0,1)\in \tau f$. \qed
}
\end{expl}

\begin{expl}
\label{expl-unfoldbeaktobeak} {\rm Let
$f\co(\RR^2,0)\to (\RR^2,0)$ be the germ given by
\[  f(x,y)=(x,y^3-x^2y). \]
Again we have $Tf=W$, where $W\subset \sE_{s,s}=\sE_{2,2}$ is the linear subspace consisting of
all $k=(k_1,k_2)$ such that the second derivative of $y\mapsto k_2(0,y)$ at $y=0$ vanishes.
The proof follows the lines of example~\ref{expl-unfoldlips}.
}
\end{expl}

\begin{expl}
\label{expl-unfoldeasyswallowtail}
{\rm Let
$f\co(\RR^2,0)\to (\RR^2,0)$ be the germ given by
\[  f(x,y)=(x,y^4+xy). \]
We want to show that $Tf$ has codimension 1 in $\sE_{s,t}=\sE_{2,2}$.
We have the standard description
\[ Tf=Jf+\tau f=\sE_s\{\,(1,y), (0,4y^3+x)\,\}+\sE_t\{\,(1,0),(0,1)\,\}. \]
A two-fold application of \cite[XV.2.1]{Martinet} proves that
$Tf+\sE_t\{\,(0,y^2)\,\}=\sE_{s,t}$.
(Follow the reasoning of example~\ref{expl-unfoldlips}.)
Hence it is enough to check that
\[  \sM_t\cdot(0,y^2) \subset Tf~. \]
As $Tf$ is an $\sE_t$-module, that reduces to checking that
\begin{align*} (0,xy^2) & \in  Tf \\
(0,y^6+xy^3) & \in  Tf~.
\end{align*}
For the first of these we write
\[ 3(0,xy^2)=4xy(1,y)-4(y^4+xy,0)+4y^4(1,y)-y^2(0,4y^3+x). \]
For the second we write
\[ 16(0,y^6+xy^3)=3x(0,4y^3+x)+4y^3(0,4y^3+x)-3x^2(0,1). \]
\qed
}
\end{expl}

\begin{expl}
\label{expl-unfoldswallowtail} {\rm
Let $g\co(\RR^2,0)\to (\RR^2,0)$ be the germ given by
\[  g(x,y)=(x,y^4+p_xy^2+xy) \]
where $x\mapsto p_x$ is a smooth function (germ) of $x$, with $p_0=0$. We shall see that
the tangent space $Tg$ has codimension $1$ in $\sE_{2,2}=\sE_{s,t}$. More precisely, we are going to
show that $g$ is left-right equivalent to the germ $f$ defined by $f(x,y)=(x,y^4+xy)$, which we investigated in
example~\ref{expl-unfoldeasyswallowtail}. Since $Tf$ has codimension 1 in $\sE_{2,2}$, it follows that
$Tg$ has codimension 1 in $\sE_{2,2}$.
\newline
For nonzero $a\in\RR$ define $g^a\co (\RR^2,0)\to (\RR^2,0)$ by
\begin{align*} g^a(x,y)=(x,y^4+a^{-2}p_{a^3x}y^2+xy). \end{align*}
Then $g^1=g$. It is easy to see that $g^a$ is left-right equivalent to $g$. Indeed, $g^a=\varphi g\psi$ where
$\psi(x,y)=(a^3x,ay)$ and $\varphi(x,y)=(a^{-3}x,a^{-4}y)$. \newline
We also define $g^0\co (\RR^2,0)\to (\RR^2,0)$ by $g^0(x,y)=(x,y^4+xy)= f(x,y)$.
With these abbreviations, the germ $G\co(\RR\times\RR^2,0) \lra (\RR\times\RR^2,0)$ defined by
\[  G(a,x,y)=(a,g^a(x,y)) \]
is smooth. (To see this, write $p_x=x\cdot u_x$ where $x\mapsto u_x$ is a smooth function. This
is possible by \cite[I.5.1]{Martinet}. Then $g^a(x,y)=(x,y^4+a\cdot u_{a^3x}y^2+xy)$, which is
clearly smooth as a function of $a$, $x$ and $y$.)
We think of it as a 1-parameter unfolding with parameter $a\in\RR$ of the germ $g^0=f$.
As $g^0$ is finitely determined, with $Tg^0$ of codimension 1 etc., we know that a miniversal
unfolding of $g^0$ is given by $F\co(\RR\times\RR^2,0) \lra (\RR\times\RR^2,0)$ where
\begin{align*}  F(b,x,y)=(b,x,y^4+by^2+xy)~. \end{align*}
By the universal property, the unfolding $G$ is isomorphic (as an unfolding of $g^0$) to the
pullback of $F$ under some map germ $\beta\co (\RR,0)\to(\RR,0)$ relating the parameter spaces.
But $\beta$ must be the zero germ. (Indeed, $g^a$ for arbitrary fixed $a$ has a serious
singularity at $0$ whereas $(x,y)\mapsto (x,y^4+by^2+xy)$ for nonzero $b$, and near the origin, has only
Whitney's folds and cusps.) Hence all $g^a$
for sufficiently small $a>0$ are left-right equivalent to $g^0=f$. But we already saw that $g^a$
for $a\ne 0$ is left-right equivalent to $g^1=g$. It follows that $g$ is left-right equivalent to $f$. \qed
}
\end{expl}

\begin{expl}
\label{expl-inflate}
{\rm
Let $f\co(\RR^2,0)\to(\RR^2,0)$ be a germ of the form
\begin{align*}   f(x,y)=(x,f_2(x,y))~. \end{align*}
Let $q\co\RR^n\to\RR$ be a nondegenerate quadratic form (a polynomial
function, homogeneous of degree $2$). Define a new
germ by
\begin{align*}  f\sha\co(\RR^{n+2},0) & \to  (\RR^2,0) \\
(z_1,\dots,z_n,x,y) & \mapsto  (x,f_2(x,y)+q(z_1,\dots,z_n)).
\end{align*}
Let $r\co \sE_{n+2,2}\to \sE_{2,2}$ be the restriction map
(restriction to the $xy$-plane). This is clearly onto. \emph{We
have}
\begin{equation}  Tf\sha= r^{-1}(Tf)~. \end{equation}
To prove this, we note first that $r(Tf\sha)\subset Tf$
and also $Tf\subset r(Tf\sha)$, from the definitions. Then it only remains to show
\[
\ker(r)\subset Tf\sha~.
\]
Indeed we shall see that $\ker(r)$ is contained in $Jf\sha$, the
subspace of $Tf\sha$ consisting of all $df\sha\cdot X$ where $X$ is
a vector field germ on $(\RR^{n+2},0)$. Suppose then that
$k=(k_1,k_2)$ is in the kernel of $r$. Let $\ell=df\sha\cdot k_1X$
where $X$ is the vector field with constant value $(0,\dots,0,1,0)$.
Then $\ell$ is in $Jf\sha\cap \ker(r)$ and $\ell_1=k_1$. Therefore
$k-\ell=(0,k_2-\ell_2)$ is in $\ker(r)$ and we only need to prove
that it is in $Jf\sha$. The function $k_2-\ell_2$ vanishes on the
$xy$-plane. Therefore, by \cite[I.5.1]{Martinet}, it can be written
in the form
\begin{align*}(z_1,\dots,z_n,x,y) & \mapsto  \sum_{i=1}^n z_i\cdot g_i(z_1,\dots,z_n,x,y).
\end{align*}
This means that $k-\ell$ can be written in the form
\begin{align*} (z_1,\dots,z_n,x,y) & \mapsto  \sum_{i=1}^n g_i(z_1,\dots,z_n,x,y)\cdot(0,z_i)~.
\end{align*}
The map $(z_1,\dots,z_n,x,y)\mapsto(0,z_i)$ is in $Jf\sha$, due to
the fact that $q$ is nondegenerate. Since $Jf\sha$ is an
$\sE_{n+2}$-submodule of $\sE_{n+2,2}$~, it follows that $k-\ell\in
Jf\sha$. \qed }
\end{expl}

\end{document}